\documentclass{amsart}

\usepackage{amsmath}
\usepackage{amsfonts}
\usepackage{amssymb}
\usepackage{amsthm,enumerate}
\usepackage[all]{xy}

\newtheorem{lem}{Lemma}[section]
\newtheorem{cor}[lem]{Corollary}
\newtheorem{prop}[lem]{Proposition}
\newtheorem{thm}[lem]{Theorem}
\newtheorem{Defn}[lem]{Definition}
\newtheorem{Ex}[lem]{Example}
\newtheorem{Question}[lem]{Question}
\newtheorem{Notation}[lem]{Notation}
\newtheorem{Property}[lem]{Property}
\newtheorem{Properties}[lem]{Properties}
\newtheorem{Subprops}{}[lem]
\newtheorem{Para}[lem]{}

\newenvironment{defn}{\begin{Defn}\rm}{\end{Defn}}
\newenvironment{ex}{\begin{Ex}\rm}{\end{Ex}}
\newenvironment{question}{\begin{Question}\rm}{\end{Question}}

\newenvironment{properties}{\begin{Properties}\rm}{\end{Properties}}
\newenvironment{subprops}{\begin{Subprops}\rm}{\end{Subprops}}
\newenvironment{para}{\begin{Para}\rm}{\end{Para}}

\newcommand{\ideal}[1]{\mathfrak{#1}}
\newcommand{\aaa}{\ideal{a}}
\newcommand{\m}{\ideal{m}}
\newcommand{\n}{\ideal{n}}
\newcommand{\p}{\ideal{p}}

\newcommand{\vf}{\varphi}

\newcommand{\rr}{\mathbb{R}}
\newcommand{\nn}{\mathbb{N}}
\newcommand{\zz}{\mathbb{Z}}

\newcommand{\x}{\mathbf{x}}

\newcommand{\HH}{\mathrm{H}}

\newcommand{\bs}{\mathrm{CI}_*}

\newcommand{\bsdim}{\bs\text{-}\mathrm{dim}}

\newcommand{\B}{\mathcal{B}}
\newcommand{\bdim}{\B\text{-}\dim}
\renewcommand{\ker}{\mathrm{Ker}}

\newcommand{\Hom}{\mathrm{Hom}}
\newcommand{\rank}{\mathrm{rank}}	
\newcommand{\cx}{\mathrm{cx}}	
	
\newcommand{\pdim}{\mathrm{pd}}
\newcommand{\pd}{\mathrm{pd}}
\newcommand{\id}{\mathrm{id}}
\newcommand{\pci}{$\mathrm{CI}_*$}
\newcommand{\pciclass}{\mathrm{CI}_*}
\newcommand{\gdim}{\mathrm{G}\text{-}\mathrm{dim}}	
\newcommand{\cidim}{\mathrm{CI}\text{-}\mathrm{dim}}

\newcommand{\pcidim}{\mathrm{CI}_*\text{-}\mathrm{dim}}	
	
\newcommand{\gldim}{\mathrm{gl}\text{-}\mathrm{dim}}
\newcommand{\glcidim}{\mathrm{glCI}\text{-}\mathrm{dim}}
\newcommand{\glpcidim}{\mathrm{glCI}_*\text{-}\mathrm{dim}}
\newcommand{\glgdim}{\mathrm{glG}\text{-}\mathrm{dim}}
\newcommand{\depth}{\mathrm{depth}} 
\newcommand{\from}{\leftarrow}
\newcommand{\im}{\text{Im}}	
\newcommand{\spec}{\mathrm{Spec}}
\newcommand{\mspec}{\mathrm{Max}}

\newcommand{\ext}{\mathrm{Ext}}	
\newcommand{\rhom}{\mathbf{R}\mathrm{Hom}}	
\newcommand{\lotimes}{\otimes^{\mathbf{L}}}

\begin{document}
\bibliographystyle{amsplain}

\author{Sean Sather-Wagstaff}
\address{Department of Mathematics, University of Illinois, 
273 Altgeld Hall, 1409 West Green Street, Urbana, IL, 61801}
\email{ssather@math.uiuc.edu}
\thanks{This research 
was conducted while the author was an NSF Mathematical Sciences 
Postdoctoral Research Fellow.}
\title{Complete intersection dimensions for complexes}
\keywords{CI-dimension, PCI-dimension, complete intersection dimension, 
homological dimension, complex}
\subjclass[2000]{13C40, 13D05, 13D25, 13D40, 14H10}  

\begin{abstract}
We extend the notions of complete intersection dimension and 
lower complete intersection dimension to the category of 
complexes with finite homology
and verify basic properties analogous to those holding for 
modules.  We also discuss the question of the behavior of complete 
intersection dimension with respect to short exact sequences.
\end{abstract}

\maketitle

\section{Introduction} \label{sec:intro}

A familiar numerical invariant of a finitely generated module over a 
Noetherian ring is its projective 
dimension.  The last few decades have seen a number of 
refinements and extensions of this.  One refinement is the notion
of Gorenstein dimension, introduced by 
Auslander and Bridger~\cite{auslander:smt}.  More recently,
Avramov, Gasharov, and Peeva~\cite{avramov:cid} defined a concept of
complete intersection dimension, Gerko~\cite{gerko:ohd} forwarded
definitions for lower complete intersection dimension and Cohen-Macaulay 
dimension, and Veliche~\cite{veliche:cmfhd} did the same for upper 
Gorenstein dimension.  The notions of complete intersection dimension 
and lower complete intersection dimension are the primary focuses of 
this paper.

These homological dimensions are well-behaved in a number of senses.
For example, when $M$ is a finite module over a Noetherian ring 
$R$ there are inequalities 
\[ \gdim_R(M)\leq\pcidim_R(M)\leq\cidim_R(M)\leq\pd_R(M); \]
if one of these dimensions is finite, then it equals those 
to its left.  When $R$ is 
local each homological dimension satisfies an 
``AB-formula'':  if one of the quantities in the displayed formula is 
finite, then it equals $\depth(R)-\depth_R(M)$.  
Furthermore, the finiteness of a
homological dimension for all finite $R$-modules
characterizes the corresponding ring-theoretic 
property of $R$ as in the theorem of Auslander, Buchsbaum, and Serre.

In another direction, the projective dimension and 
Gorenstein 
dimension have been extended to complexes of $R$-modules.  
The projective dimension was systematically developed by
Foxby~\cite{foxby:hacr,foxby:ibc,foxby:bcofm},
and the G-dimension by 
Yassemi~\cite{yassemi:gd} and Christensen~\cite{christensen:gd}.
The purpose of this paper is to give a similar extension of 
complete intersection dimension 
and lower complete intersection dimension and
verify basic properties that one 
expects to carry over from the situation for modules.
This is done in Sections~\ref{sec:cidim} and~\ref{sec:pcidim}.  
Also, we prove 
stability results, Theorems~\ref{thm:cidim-stability}
and~\ref{thm:pcidim-stability}, that are particular to complexes.

One difficulty with the complete intersection dimension is that we do
not know whether it is well-behaved with respect to short 
exact sequences; Section~\ref{sec:es} is devoted to this issue.  
Section~\ref{sec:gldim} consists of a brief discussion of ``global'' 
homological dimensions, which can be introduced from the homological 
dimensions under consideration, like the global 
dimension of Cartan and Eilenberg~\cite{cartan:ha}.  
Section~\ref{sec:back} is home to a brief catalogue of background 
material used in the other sections.

Acknowledgements.  I am grateful to L.~Avramov, A.~Frankild, and 
S.~Iyengar for helpful discussions about this material and its 
presentation.

\section{Background} \label{sec:back}

This section is mostly a summary of standard notions from 
hyperhomological algebra;  the interested reader is 
directed to~\cite{foxby:hacr} for a detailed account.
We also include a couple of results that will be 
important in the sections that follow. 

\textbf{Throughout this paper, all rings are commutative and 
Noetherian.}  

A complex of modules over a ring $R$ is a 
sequence of $R$-module homomorphisms
\[ X=\cdots\xrightarrow{\partial^X_{i+1}}X_i
\xrightarrow{\partial^X_{i}}X_{i-1}\xrightarrow{\partial^X_{i-1}}
\cdots \]
such that $\partial^X_i\partial^X_{i+1}=0$ for every integer $i$.  
When $M$ 
is an $R$-module, identify $M$ with the complex 
$\cdots\to 0\to M\to 0\to\cdots$ concentrated in degree 0.  

\begin{para} \label{para:back-1}
A complex $X$ is \emph{bounded below} (resp., \emph{bounded}) if 
$X_i=0$ for all $i\ll 0$ (resp., for all $|i|\gg 0$);  it is 
\emph{degreewise finite} if each $X_i$ is a finite $R$-module;  and 
it is \emph{finite} if it is bounded and degreewise finite.
Next, $X$ is 
\emph{homologically bounded below} (resp., \emph{homologically bounded}) 
if the homology complex $\HH(X)$ is bounded below (resp., bounded);  
it is \emph{homologically degreewise finite} (resp., 
\emph{homologically finite}) if $\HH(X)$ is degreewise finite (resp., 
bounded and degreewise finite).  The 
\emph{supremum} and \emph{infimum} 
of $X$ are given by the following 
formulas:
\[ \sup(X)=\sup\{i\in\zz\mid\HH_i(X)\neq 0\}\quad\text{and}\quad
\inf(X)=\inf\{i\in\zz\mid\HH_i(X)\neq 0\}. \]
Given an integer $n$, the \emph{$n$th suspension} 
of $X$ is the complex $\Sigma^nX$ with $(\Sigma^nX)_m=X_{m-n}$ and 
differential $\partial_m^{\Sigma^nX}=(-1)^n\partial^X_{m-n}$ for each 
$m$.  The 
kernel and cokernel of $\partial^X_{n}$ are denoted
$Z^X_{n}$ and $C^X_{n-1}$, respectively.  For any $R$-module $M$, 
one has $C^{X\otimes_R M}_n\cong C^X_n\otimes_R M$, 
by the right-exactness of $\text{-}\otimes_R M$.
The \emph{$n$th 
soft left- and right-truncations} of $X$ are the complexes
\begin{align*}
\tau_{\leq n}(X)&=\cdots\to 0\to 
C^X_n\xrightarrow{\overline{\partial^X_n}} 
X_{n-1}\xrightarrow{\partial^X_{n-1}}\cdots \\
\tau_{\geq n}(X) &=\cdots\xrightarrow{\partial^X_{n+2}} 
X_{n+1}\xrightarrow{\partial^X_{n+1}} Z^X_n\to 0\to\cdots \\
\end{align*}
respectively,
where $\overline{\partial^X_n}$ is the map induced by $\partial^X_n$.  
The \emph{$n$th hard left- and right-truncations} are the 
complexes
\begin{align*}
X_{\leq n}&=\cdots\to 0\to X_n\xrightarrow{\partial^X_n} 
    X_{n-1}\xrightarrow{\partial^X_{n-1}}\cdots \\
X_{\geq n}&=\cdots\xrightarrow{\partial^X_{n+2}} X_{n+1}
    \xrightarrow{\partial^X_{n+1}} X_n\to 0\to\cdots.
\end{align*}
\end{para}

It is worth noting explicitly that we do not use the machinery of derived 
categories in this paper.  This is for two reasons:
we are interested in how the invariants we define behave with 
respect to short exact sequences, and we use kernels and 
cokernels of morphisms in our arguments.  Instead, we work within the 
category of complexes of modules.

\begin{para} \label{para:back-2}
Let $X,Y$ be complexes of $R$-modules.  A \emph{morphism} 
$\sigma\colon X\to Y$ is a collection of $R$-module homomorphisms 
$\sigma_i\colon X_i\to Y_i$ such that 
$\partial_i^Y\sigma_i=\sigma_{i-1}\partial_i^X$ for each integer $i$.
A \emph{quasiisomorphism} is a 
morphism $\alpha\colon X\to Y$ such that the
map induced on  homology 
$\HH(\alpha)\colon\HH(X)\to\HH(Y)$ is an isomorphism;  this is 
signified by $\alpha\colon X\xrightarrow{\simeq} Y$.  
More generally, $X$ and $Y$ are \emph{quasiisomorphic}, denoted 
$X\simeq Y$ if there is a finite sequence of quasiisomorphisms
\[ X\xleftarrow{\simeq}X^1\xrightarrow{\simeq}X^2\xleftarrow{\simeq} 
\cdots\xrightarrow{\simeq}Y.\]
If $m\leq\inf(X)$ and $n\geq\sup X$, then the natural maps 
$X\to\tau_{\leq n}(X)$, $\tau_{\geq m}(X)\to X$, and 
$X_{\geq n}\to\Sigma^nC^X_n$ are quasiisomorphisms.  
Thus, if $s=\sup(X)$, then $C^X_s\neq 0$.
\end{para}

The homological dimensions studied in this work are descendants of the 
projective dimension.

\begin{para} \label{para:back-3}
A \emph{projective} (resp., \emph{free}) \emph{resolution} 
of a homologically bounded below complex $X$ is a bounded below 
complex $P\simeq X$ of projective (resp., free) $R$-modules.  
If $X$ is homologically both degreewise finite and bounded below, 
then it possesses a degreewise finite free resolution; 
see~\cite[(2.6.L)]{foxby:hacr} or apply~\cite[3.1.6]{roberts:maccila} 
to the truncation $\tau_{\geq m}(X)\simeq X$ for $m=\inf(X)$.  
By~\cite[(1.2.P, 1.4.P)]{avramov:hdouc}, if $P\simeq X$ is a 
projective resolution, then there exists a quasiisomorphism 
$P\xrightarrow{\simeq} X$.  

The \emph{projective dimension of $X$} is
\[ \pd_R(X)=\inf\{\sup\{n\mid P_i\neq 0\}\mid\text{$P$ is a projective 
resolution of $X$}\}.\]
Thus, if $\pd_R(X)$ is finite, then $X$ is 
homologically both bounded and nonzero.
Injective resolutions and the injective dimension $\id_R(X)$ are 
defined dually.
\end{para}

Given a morphism of complexes $X\to Y$ it can be useful to be able 
to enlarge $X$ to construct a surjective morphism with the same
morphism induced on homology.  
The next fact~\cite[(8.4.4,5)]{avramov:dgha} allows us to do so.  
See~\ref{lem:back-1} and~\ref{prop:ex-seq-resol} for applications.

\begin{para} \label{para:back-n}
Given a bounded below degreewise finite complex of $R$-modules $X$, 
there exists a bounded below degreewise finite complex of free 
$R$-modules $G$ with $\HH(G)=0$ and a morphism 
$\epsilon\colon G\to X$ such that each $\epsilon_i$ is surjective.
\end{para}

The following is a version of the 
existence of ``strict semifree resolutions''.

\begin{lem} \label{lem:back-1}
Let $R$ be a ring and $X$ a complex of $R$-modules that is bounded 
below and degreewise finite.  There exists 
a degreewise finite free resolution $\sigma\colon P\xrightarrow{\simeq} X$ such 
that each $\sigma_i\colon P_i\to X_i$ is surjective.
\end{lem}

\begin{proof}
By~\cite[3.1.6]{roberts:maccila} take a
degreewise finite free resolution 
$\alpha\colon F\xrightarrow{\simeq} X$.  
Fix a complex $G$ and morphism $\epsilon\colon G\to X$ as
in~\ref{para:back-n}.
The complex $P=F\oplus G$ and morphism
$\sigma\colon P\to X$ given by 
$\sigma_i(f,g)=\alpha_i(f)+\epsilon_i(g)$ satisfy
the conclusions.
\end{proof}

Given a short exact sequence of complexes, it is well-known that 
there exists a short exact sequence on the level of projective 
resolutions~\cite[($6.10^{\circ}$)]{iversen:cs}.  
It is helpful to know when the projective resolutions can be chosen 
to be degreewise finite.

\begin{prop} \label{prop:ex-seq-resol}
Let $R$ be a ring and 
$0\to X\xrightarrow{\eta} Y\xrightarrow{\nu} Z\to 0$
an exact sequence 
of complexes of $R$-modules that are homologically both degreewise 
finite and bounded below.  There exists
a commutative diagram of complexes with exact rows 
\[ \xymatrix{0\ar[r] & T\ar[r]\ar[d]_{\simeq}^{\psi} & 
U\ar[r]\ar[d]_{\simeq}^{\lambda} & V\ar[r]\ar[d]_{\simeq}^{\alpha} & 0 \\
0\ar[r] & X\ar[r]^{\eta} & Y\ar[r]^{\nu} & X\ar[r] & 0
} \]
where each vertical map is a
degreewise finite $R$-projective resolution.
\end{prop}

\begin{proof}
Let $\alpha\colon V\xrightarrow{\simeq} Z$ and 
$\gamma\colon F\xrightarrow{\simeq} Y$
be degreewise finite $R$-free resolutions.  
There exists a 
morphism $\sigma\colon F\to V$ such that $\nu\gamma=\alpha\sigma$, 
by~\cite[(1.1.P.1),(1.2.P)]{avramov:hdouc}.  
Since $V$ is bounded below and degreewise finite, fix a complex 
$G$ and morphism $\epsilon\colon G\to V$ as in~\ref{para:back-n}.
By~\cite[(9.8.3.2'),(9.7.1)]{avramov:dgha} there exists a 
morphism $\rho\colon G\to Y$ such that $\nu\rho=\alpha\epsilon$.  

Let $U=F\oplus G$ and define morphisms $\lambda\colon U\to Y$ and 
$\theta\colon U\to V$ by the formulas
\[ \lambda_i(f,g)=\gamma_i(f)+\rho_i(g)\qquad\text{and}\qquad
\theta_i(f,g)=\sigma_i(f)+\epsilon_i(g).\]
It is straightforward to check that $\alpha\theta=\nu\lambda$.
Furthermore, $\lambda\colon U\to Y$ is a degreewise finite $R$-free 
resolution and each
$\theta_i$ is surjective.  Set $T=\ker(\theta)$ with 
$\iota\colon T\to U$ the natural inclusion and $\psi\colon T\to X$ 
the morphism induced  by $\lambda$.  Since each sequence 
$0\to T_i\to U_i\to V_i\to 0$ is exact with $U_i,V_i$ projective,
each $T_i$ is projective.  Thus, we have a commutative 
diagram of the desired form.
The 5-lemma applied to the long exact 
sequences in homology shows that $\psi$ is a quasiisomorphism.
\end{proof}

Derived Hom and tensor product are ubiquitous tools in the study of 
complexes.

\begin{para} \label{para:back-5}
Given complexes of $R$-modules $X,Y$ 
with $X$ homologically bounded below, then
$X\lotimes_R Y$ and $\rhom_R(X,Y)$ denote
the complexes 
$P\otimes_R Y$ and $\Hom_R(P,Y)$, respectively, where $P\simeq X$ is a 
projective resolution.  These complexes are 
only well-defined up to quasiisomorphism, but this is enough for our 
applications.
\end{para}

The G-dimension
comes to bear directly and indirectly on the study of 
complete intersection dimension.  A nice treatment can be 
found in~\cite{christensen:gd}.

\begin{para} \label{para:back-4}
For a ring $R$, let $(\text{-})^*=\Hom_R(\text{-},R)$.
A finite $R$-module $M$ is \emph{totally reflexive over $R$} if 
$M$ is reflexive and $\ext^i_R(M,R)=0=\ext^i_R(M^*,R)$ for all $i>0$.  
Each finitely generated projective 
$R$-module is totally reflexive over $R$.
A \emph{G-resolution} of a 
complex $X$ is a bounded below complex $G\simeq X$, 
such that each $G_i$ is totally reflexive over $R$.  The 
\emph{G-dimension} of $X$ is
\[ \gdim_R(X)=\inf\{\sup\{n\mid G_i\neq 0\}\mid\text{$G$ is a 
G-resolution of $X$}\}.\]
By~\cite[(2.3.8)]{christensen:gd}, if $\gdim_R(X)<\infty$, then one has
\[ \gdim_R(X)=-\inf(\rhom_R(X,R)). \]
\end{para}

The depth of a finite module over over a local ring is a familiar 
invariant.  Our definition of depth for complexes is taken from 
Iyengar~\cite{iyengar:dfcait}.

\begin{para} \label{para:back-6}
Let $R$ be a local ring and $K$ the Koszul complex over $R$ on a 
sequence of generators of length $n$ for the maximal ideal $\m$
of $R$.  For a complex of $R$-modules $X$, the \emph{depth of $X$} is
\[ \depth_R(X)=n-\sup(X\otimes_R K).\]
This is independent of the sequence of generators for $\m$.
\end{para}

Complexes of finite projective dimension have finite G-dimension, and 
the finiteness of either of these implies
an AB-formula, where ``AB'' stands for Auslander-Buchsbaum 
and Auslander-Bridger, c.f.~\cite[(2.3.10,13)]{christensen:gd}.

\begin{para} \label{para:AB}
For a homologically finite complex $X$ over a ring $R$, one has an
inequality
\[ \gdim_R(X)\leq\pdim_R(X) \]
with equality when $\pdim_R(X)<\infty$.  If $R$ is local and 
$\gdim_R(X)<\infty$, then
\[ \gdim_R(X)=\depth(R)-\depth_R(X).\]
\end{para}

The Betti numbers of a complex over a local ring are of particular 
interest in connection with the complete intersection dimension.

\begin{para} \label{para:back-7}
Let $(R,\m,k)$ be a local ring and $X$ a homologically bounded below and
degreewise finite complex of $R$-modules.  
By~\cite[(2.2.4)]{roberts:hiomocr}, $X$ has a 
\emph{minimal} free resolution, that is, a degreewise finite free 
resolution $F\simeq X$ such that $\partial^F(F)\subseteq\m F$.  As 
is the case with modules, minimal free resolutions are unique up to 
isomorphism.  The \emph{$n$th Betti number} of $X$ is
\[ \beta_n^R(X):=\rank_R(F_n)=\rank_k\HH_n(X\lotimes_R k). \]
The \emph{Poincar\'{e} series} of $X$ is the formal Laurent series
\[ P_X^R(t)=\sum_n\beta_n^R(X)t^n.\]  
The \emph{complexity} of $X$, defined by the formula
\[ \cx_R(X)=\inf\{c\in\nn\mid \text{there exists $\alpha\in\rr$ such 
that $\beta^R_n(X)\leq \alpha n^{c-1}$ for $n\gg 0$}\} \]
is a measure of the asymptotic size of the minimal free resolution of 
$X$.  For instance, $\cx_R(X)=0$ if and only if $\pd_R(X)<\infty$.
\end{para}

The behavior ``at infinity'' of the sequence of 
Betti numbers of a complex is 
almost identical to that of the syzygy modules of the complex.

\begin{para} \label{para:back-z}
Let $X$ be a homologically finite complex of modules over a local 
ring $R$, and fix a degreewise finite $R$-free resolution $P\simeq X$.
For $n\geq \sup(X)$, it is straightforward to show that the 
Poincar\'{e} series of $X$ and $C^P_n$ are related by the formula
\[ P_X^R(t)=t^nP_{C^P_n}^R(t)+t^{\inf(X)}f(t) \]
for some polynomial 
$f(t)\in\zz[t]$.  In particular, it follows that $\cx_R(X)=\cx_R(C^P_n)$.
\end{para}

Certain accounting principles~\cite[(11.11)]{foxby:hacr}
are handy for tracking the behavior 
of complexity under derived tensor product.

\begin{para} \label{para:back-cx-add}
Let $R$ be a local ring with homologically finite 
complexes $X,Y$.  There is an equality of Poincar\'e series
\[ P^R_{X\lotimes_R Y}(t)=P^R_X(t)P^R_Y(t).\]
It follows that, if $\HH(Y)\neq 0$, then 
\[ \cx_R(X)\leq\cx_R(X\lotimes_R Y)\leq\cx_R(X)+\cx_R(Y).\]
In particular, if $X$ and $Y$ have finite complexity, then so has 
$X\lotimes_R Y$.
\end{para}

\section{Complete Intersection Dimension for Complexes} \label{sec:cidim}

In this section, we introduce the notion of CI-dimension for 
homologically finite complexes and  
verify a number of properties which the CI-dimension
for modules leads us to expect.  For a nonzero finite module, 
considered as a complex concentrated in degree 0, the 
definition is the same as that given in~\cite{avramov:cid}.

\begin{defn}  \label{def:cidim-1}
Let $R$ be a ring and $X$ a homologically finite complex of 
$R$-modules.  
When $R$ is local, a
\emph{(codimension $c$) quasi-deformation} of $R$ is a diagram 
of local homomorphisms $R\to R'\from Q$ such that 
the first map is flat and the second 
map is surjective with kernel generated by a 
$Q$-sequence (of length $c$).  In this situation, 
let $X'$ denote the complex $X\otimes_R R'$.
The \emph{CI-dimension of $X$} is
\[ \cidim_R(X)=\inf\{\pdim_Q(X')-\pdim_Q(R')\mid
\text{$R\to R'\from Q$ is a quasi-deformation} \}. \]
When $R$ is not necessarily local the 
\emph{CI-dimension of $X$} is
\[ \cidim_R(X)=\sup\{\cidim_{R_{\m}}(X_{\m})\mid\m\in\mspec(R) \} \]
where $\mspec(R)$ is the set of all maximal ideals of $R$.  
\end{defn}

Certain facts are immediate from the definition.  

\begin{properties} \label{props:cidim}
Fix a ring $R$ and a 
homologically finite complex of $R$-modules $X$.
\begin{subprops} \label{subprop:cidim-1} 
$\cidim_R(X)\in\{-\infty\}\cup\zz\cup\{\infty\}$.
\end{subprops}  
\begin{subprops} \label{subprop:cidim-2}
$\cidim_R(X)=-\infty$ if and only if $X\simeq 0$.
\end{subprops}  
\begin{subprops} \label{subprop:cidim-5}
If $X\simeq Y$, then $\cidim_R(X)=\cidim_R(Y)$. 
\end{subprops}  
\begin{subprops} \label{subprop:cidim-4}
Each integer $n$ yields $\cidim_R(\Sigma^nX)=\cidim_R(X)+n$.
\end{subprops}  
\end{properties}

The CI-dimension for complexes
fits into a hierarchy of homological dimensions 
like that for modules.  Also, over a 
local ring, an AB-formula is satisfied.  This is the analogue 
of~\cite[(1.4)]{avramov:cid} for complexes;  since the proof is 
identical, we omit it here.

\begin{prop} \label{prop:cidim-1}
Let $R$ be a ring and $X$ a homologically finite complex of $R$-modules.  
There are inequalities
\[ \gdim_R(X)\leq\cidim_R(X)\leq\pdim_R(X); \]
when one of these dimensions is finite it is equal to those 
on its left.  In particular, 
$\sup(X)\leq\cidim_R(X)$.
If $R$ is local and $\cidim_R(X)<\infty$, then 
\begin{xxalignat}{3}
&{\hphantom{\square}}& \cidim_R(X) &=\depth(R)-\depth_R(X).
&&\square
\end{xxalignat}
\end{prop}

Like the G-dimension and projective dimension,
CI-dimension is well-behaved with respect to localization.  
Again, the proof is identical to that of the corresponding result for 
modules~\cite[(1.6)]{avramov:cid}

\begin{prop} \label{prop:cidim-2}
Let $R$ be a ring and $X$ a homologically finite complex of 
$R$-modules.  For every 
multiplicative subset $S\subset R$, there is an inequality
\[ \cidim_{S^{-1}R}(S^{-1}X)\leq\cidim_R(X).\]
Furthermore, 
\begin{xxalignat}{3}
&{\hphantom{\square}}& \cidim_R(X)&
=\sup\{\cidim_{R_{\p}}(X_{\p})\mid\p\in\spec(R)\}.
&&\square
\end{xxalignat}
\end{prop}

The following proposition is the expected analogue of
the Avramov-Gasharov-Peeva characterization of local complete 
intersection rings~\cite[(1.3)]{avramov:cid}.  
Recall that a ring $R$ is ``locally 
a complete intersection'' if, for every maximal ideal $\m$ of $R$, 
the localization $R_{\m}$ is a complete intersection.

\begin{prop} \label{prop:cidim-3}
For a ring $R$ with $\dim(R)<\infty$,  
the following are equivalent.
\begin{enumerate}[(a)]
\item  \label{item:cidim-1}
$R$ is locally a complete intersection.
\item \label{item:cidim-2}
Each homologically finite complex of $R$-modules $X$ satisfies
\[ \cidim_R(X)\leq\dim(R)+\sup(X).\]
\item  \label{item:cidim-3}
Each maximal ideal $\m\subset R$ satisfies
$\cidim_R(R/\m)<\infty$.
\end{enumerate}
\end{prop}

\begin{proof}
``(\ref{item:cidim-1})$\implies$(\ref{item:cidim-2})''.  
Let $X$ be a homologically finite complex of 
$R$-modules.  
For each maximal ideal $\m$ of $R$, one has
$\cidim_{R_{\m}}(X_{\m})<\infty$;  the proof is identical to that 
of~\cite[(1.3)]{avramov:cid}.  Furthermore, 
\[ \cidim_{R_{\m}}(X_{\m})=\depth(R_{\m})-\depth_{R_{\m}}(X_{\m})
\leq\dim(R_{\m})+\sup(X_{\m}) \]
where the equality is by~\ref{prop:cidim-1} and the 
inequality is by~\cite[(2.7)]{foxby:daafuc}.  It follows that
\begin{align*}
\cidim_R(X)&=\sup\{\cidim_{R_{\m}}(X_{\m})\mid\m\in\mspec(R)\} \\
&\leq\sup\{\dim(R_{\m})+\sup(X_{\m})\mid\m\in\mspec(R)\} \\
&\leq\sup\{\dim(R_{\m})\mid\m\in\mspec(R)\}+
\sup\{\sup(X_{\m})\mid\m\in\mspec(R)\}\\
&=\dim(R)+\sup(X).
\end{align*}

``(\ref{item:cidim-2})$\implies$(\ref{item:cidim-3})'' is trivial.

``(\ref{item:cidim-3})$\implies$(\ref{item:cidim-1})''.  
By definition
$\cidim_{R_{\m}}(R_{\m}/\m R_{\m})=\cidim_R(R/\m)<\infty$, and so
$R_{\m}$ is a complete intersection by~\cite[(1.3)]{avramov:cid}.
\end{proof}

The next result is the main tool used to understand the 
relation between the CI-dimension of a complex $X$ and 
that of its syzygy 
modules.

\begin{lem} \label{lem:cidim-1}
Let $R$ be a ring and 
$0\to X^1\to X^2\to X^3\to 0$ an exact sequence of 
homologically finite complexes of $R$-modules.
For integers $i,j,k$ such that $\{i,j,k\}=\{1,2,3\}$, there is an 
inequality
\[ \cidim_R(X^k)\leq\max\{\pdim_R(X^i),\cidim_R(X^j)\}+1. \]
In particular, if $\pdim_R(X^i)$ and $\cidim_R(X^j)$ are
finite, then $\cidim_R(X^k)<\infty$.
\end{lem}

\begin{proof}
Assume that $\pdim_R(X^i),\cidim_R(X^j)<\infty$
and pass to $R_{\m}$ to assume that
$R$ is local.  Let $R\to R'\from Q$ be a 
codimension $c$ quasi-deformation 
such that $\pdim_Q((X^j)')<\infty$.  
It is straightforward to show that
\[ \pdim_Q((X^k)')
    \leq\max\{\pdim_Q((X^i)'),\pdim_Q((X^j)')\}+1<\infty. \]
The desired conclusion now follows from 
the equalities $\pdim_R(X^i)=\cidim_R(X^i)$ and
$\cidim_R(X^m)=\pdim_Q((X^m)')-c$ for $m=1,2,3$
\end{proof}

Given an exact sequence as in the lemma, it is not known whether one 
can replace $\pdim_R(X^i)$ with $\cidim_R(X^i)$, even when each complex is 
a module concentrated in degree 0.  This issue is discussed further in 
Section~\ref{sec:es}.

As is the case for modules~\cite[(1.9)]{avramov:cid},
one can compute the CI-dimension of a complex from
that of its syzygies.

\begin{prop} \label{prop:cidim-4}
Let $X$ be a homologically finite complex of $R$-modules.
Fix a degreewise finite $R$-projective resolution 
$P\simeq X$ and an integer $n\geq\sup(X)$.
\begin{enumerate}[(i)]
\item  If $C^P_n= 0$, then $\cidim_R(X)=\pdim_R(X)<n<\infty$.
\item  If $C^P_n\neq 0$, then $\cidim_R(C^P_n)=\max\{0,\cidim_R(X)-n\}$.
\end{enumerate}
\end{prop}

\begin{proof}
Consider the exact sequence of complexes
\begin{equation} 
0\to P_{\leq n-1}\to P\to P_{\geq n}\to 0 \tag{$\star$} \label{seq}
\end{equation}
and recall that $P_{\geq n}\simeq\Sigma^n C^P_n$.
If $C^P_n=0$, then the morphism
$P_{\leq n-1}\to P$ is a quasiisomorphism, and it follows that
$\pdim_R(X)=\pdim_R(P)=\pdim_R(P_{\leq n-1})< n$.

If $C^P_n\neq 0$, then 
$\cidim_R(X)<\infty$ if and only if  $\cidim_R(P_{\geq n})<\infty$
by Lemma~\ref{lem:cidim-1}.  Since 
$\cidim_R(P_{\geq n})=\cidim_R(C^P_n)+n$, the formula holds when 
$\cidim_R(X)=\infty$, so assume that $\cidim_R(X)<\infty$.
The CI-dimensions of the complexes in (\ref{seq})
agree with their G-dimensions.  An analysis of the long exact 
sequence on homology associated to the  exact sequence
$\rhom(\text{(\ref{seq})},R)$ shows that
\begin{align*}
\gdim_R(C^P_n)+n 
  & = \gdim_R(P_{\geq n}) \\
  & = -\inf(\Hom_R(P_{\geq n},R)) \\
  & = -\min\{-n,\inf(\Hom_R(P,R))\} \\
  & = \max\{n,\gdim_R(P)\} 
\end{align*}
and the result now follows from Proposition~\ref{prop:cidim-1}.
\end{proof}

\begin{cor} \label{cor:cidim-4-0}
For a homologically finite complex of $R$-modules $X$
the following conditions are equivalent.
\begin{enumerate}[(a)]
\item $\cidim_R(X)< \infty$.
\item Each degreewise finite $R$-projective resolution $P\simeq X$
and each $n\geq\sup(X)$ yield $\cidim_R(C^P_n)< \infty$.  
\item Some degreewise finite $R$-projective resolution $P\simeq X$ 
and some $n\geq\sup(X)$ yield $\cidim_R(C^P_n)< \infty$.  
\qed
\end{enumerate}
\end{cor}

As a corollary, one sees that a complex of finite CI-dimension has what 
might be termed a ``finite CI-resolution''.  The converse of this 
property is related to the behavior of CI-dimension over short exact 
sequences;  see Theorem~\ref{thm:es}.

\begin{cor} \label{cor:es}
If $\cidim_R(X)$ is 
finite, then there exists 
a finite complex of $R$-modules $Y\simeq X$ such that each 
nonzero $Y_i$ has CI-dimension 0.
\end{cor}

\begin{proof}
Let $n=\cidim_R(X)$ and fix a degreewise finite projective resolution 
$P\simeq X$.  Consider the soft truncation 
$\tau_{\leq n}(P)\simeq X$.
Then $\tau_{\leq n}(P)_i=0$ for each $i>n$ and $\tau_{\leq n}(P)_i$ is a
finitely generated projective for each 
$i\neq n$.  Proposition~\ref{prop:cidim-4} 
implies that $\tau_{\leq n}(P)_n\cong 
C_n^P$ has CI-dimension 0 so that $\tau_{\leq n}(P)$ has the desired form.
\end{proof}

We now
use Proposition~\ref{prop:cidim-4} 
to deduce facts about 
CI-dimension for complexes directly from the corresponding 
facts for 
modules~\cite[(1.12,13),(4.10),(5.3,6)]{avramov:cid}.  
It is worth noting that the results on complexity can be proved
using cohomological operators as 
in~\cite{avramov:codbad} and~\cite{yoshino}.

\begin{cor} \label{cor:cidim-2}
Let $X$ be a homologically finite complex of 
$R$-modules.
\begin{enumerate}[(i)]
\item  \label{item:cidim-7}
For a faithfully flat ring homomorphism $R\to S$ there is an 
inequality
\[ \cidim_R(X)\leq\cidim_S(X\otimes_R S) \]
with equality when $\cidim_S(X\otimes_R S)<\infty$.
\item \label{item:cidim-8}
Let $\pi\colon Q\to R$ be a surjective ring homomorphism with kernel 
generated by a $Q$-regular sequence $\x=x_1,\ldots,x_c$.  There is an
inequality
\[ \cidim_R(X)\leq\cidim_Q(X)-c \]
with equality when $\cidim_Q(X)<\infty$.
\item  \label{item:cidim-9}
Let $\aaa\subset R$ be an ideal, $R^*$ the $\aaa$-adic 
completion, and $X^*=X\otimes_R R^*$.  There is an
inequality
\[ \cidim_{R^*}(X^*)\leq\cidim_R(X) \]
with equality when $\aaa$ is contained in the Jacobson radical of $R$.
\item  \label{item:cidim-10}
If $R$ is local and $\cidim_R(X)$ finite, then the Poincar\'{e} series 
$P^R_X(t)$ is a rational function in $\zz(t)$, and $\cx_R(X)$ is equal to 
the order of the pole at $t=1$ of $P^R_X(t)$;  in 
particular, $\cx_R(X)<\infty$.
\item  \label{item:cidim-11}
If $R$ is local and $\cidim_R(X)<\infty$, then 
$\cx_R(X)\leq\mathrm{edim}(R)-\depth(R)$, and the inequality is strict 
unless $R$ is a complete intersection.
\end{enumerate}
\end{cor}

\begin{proof}
(\ref{item:cidim-7})  If $P$ is a degreewise finite
$R$-free resolution of $X$, then $P\otimes_R S$ is a
degreewise finite
$S$-free resolution of $X\otimes_R S$, and 
$C_n^{(P\otimes_R S)}=C_n^P\otimes_R S$
for each integer $n$.  By~\cite[(1.13.1)]{avramov:cid} 
$\cidim_R(C^P_n)\leq\cidim_S(C^P_n\otimes_R S)$ with equality when 
$\cidim_S(C^P_n\otimes_R S)<\infty$.  
Applying~\ref{prop:cidim-4} 
with $n=\sup(X)$ 
implies the desired result.

(\ref{item:cidim-8})  Assume that $\cidim_Q(X)<\infty$.  
For a maximal ideal $\n$ of $Q$ not containing $\x$, one has 
$X_{\n}=0$.  Thus, one reduces to the case where $Q$ and $R$ are 
local.  In this case, apply~\ref{prop:cidim-4} 
with $n=\sup(X)$ 
and~\cite[(1.12.3)]{avramov:cid} as in (\ref{item:cidim-7}), 
to deduce the result.

(\ref{item:cidim-9})  This is proved similarly to 
(\ref{item:cidim-7}), using~\cite[(1.13.2)]{avramov:cid}.

(\ref{item:cidim-10})
Let $P$ be a minimal free resolution of $X$ and fix an integer 
$n\geq\sup(X)$.  By~\ref{prop:cidim-4}, 
$\cidim_R(C^P_n)<\infty$.  By~\cite[(4.10)]{avramov:cid} 
and~\cite[(11.1)]{atiyah:ica}, the Poincar\'{e} series $P_{C^P_n}^R(t)$ is 
in $\zz(t)$.  By~\cite[(5.3)]{avramov:cid}, the order of 
the pole of $P_{C^P_n}^R(t)$ at $t=1$ is exactly $\cx_R(C^P_n)$.  
As noted in~\ref{para:back-z}, one has
$P_X^R(t)=t^nP_{C^P_n}^R(t)+t^{\inf(X)}f(t)$ for some polynomial 
$f(t)\in\zz[t]$.  In particular, $P^R_X(t)\in\zz(t)$, 
the orders of the poles at $t=1$ of $P_X^R(t)$ and $P_{C^P_n}^R(t)$
are equal, and $\cx_R(X)=\cx_R(C^P_n)$.

(\ref{item:cidim-11})  
Use the equality $\cx_R(X)=\cx_R(C_n^P)$ and~\cite[(5.6)]{avramov:cid}.
\end{proof}

The final result of this section parallels stability results of 
Yassemi~\cite[(2.14,15)]{yassemi:gd} for G-dimension and their 
generalizations~\cite[(5.1,7--9)]{iyengar:gdolh}.  It is particular 
to complexes because, when $M$ and $N$ are finite modules with 
$\pdim_R(N)$ finite, the complexes $M\lotimes_R N$ and 
$\rhom_R(N,M)$ are generally not
concentrated homologically in any single degree.
Also, it is easy to construct examples showing that the hypothesis 
``$\pd_R(P)$ is finite'' 
is necessary:  even for two finite modules $M,N$ over a 
local complete intersection, $M\lotimes N$ need not be homologically 
bounded and therefore need not have finite CI-dimension.

\begin{thm} \label{thm:cidim-stability}
Let $R$ be a ring and $X,P$ homologically finite complexes of 
$R$-modules.
If $\pd_R(P)$ is finite 
then 
\begin{align*}
\cidim_R(X\lotimes_R P)&=\cidim_R(X)+\cidim_R(P) & \text{and} \\
\cidim_R(\rhom_R(P,X))&=\cidim_R(X)-\inf(P).
\end{align*}
In particular, the CI-dimensions of the complexes
$X$, $X\lotimes_R P$, and 
$\rhom_R(P,X)$ are simultaneously finite.
\end{thm}

\begin{proof}
By~\cite[(2.14,15)]{yassemi:gd}, it suffices to show that 
$\cidim_R(X)$, $\cidim_R(X\lotimes_R P)$, and 
$\cidim_R(\rhom_R(P,X))$ are simultaneously finite.  
Furthermore, it suffices to consider the case where $R$ is local and 
$\HH(X)\neq 0$.  It is straightforward to show that
$\HH(X\lotimes_R P)$ and $\HH(\rhom_R(P,X))$ are both nonzero.

For any quasi-deformation $R\to R'\from Q$ one has
$(X'\lotimes_{R'}P')\simeq(X\lotimes_R P)'$.  Since 
$\pd_{R'}(P')=\pd_R(P)<\infty$, it follows 
from~\cite[(5.8)]{iyengar:gdolh} that 
$\pd_Q((X\lotimes_R P)')=\pd_Q(X')+\pd_{R'}(P')$.  In particular, 
$\pd_Q((X\lotimes_R P)')$ and $\pd_Q(X')$ are simultaneously finite, 
and thus the same is true of 
$\cidim_R(X\lotimes_R P)$ and $\cidim_R(X)$.

The tensor-evaluation morphism
$X\lotimes_R\rhom_R(P,R)\to\rhom_R(P,X)$ is a quasiisomorphism, 
because $\pd_R(P)<\infty$ and $\HH(P)$ is finite.  
Since $\rhom_R(P,R)$ is homologically finite
and $\pd_R(\rhom_R(P,R))<\infty$, the 
last paragraph implies that
$\cidim_R(\rhom_R(P,X))$ is finite if and only if $\cidim_R(X)$ is 
finite.
\end{proof}

\section{Exact Sequences} \label{sec:es}

In this section, we discuss the behavior of 
CI-dimension with respect to exact sequences.  
The primary question is the following.

\begin{question} \label{quest:es}
Let $R$ be a ring and $0\to X^1\to X^2\to X^3\to 0$ an exact sequence 
of homologically finite complexes of $R$-modules.  For integers $i,j,k$ 
such that $\{i,j,k\}=\{1,2,3\}$,
if $\cidim_R(X^i),\cidim_R(X^j)<\infty$, must it be that 
$\cidim_R(X^k)<\infty$?
\end{question}

For a ring $R$, if the answer to Question~\ref{quest:es}
is always ``yes'', the ring
$R$ is said to \textit{satisfy the exact sequence property (ES)}.  If 
the answer is always ``yes'' for exact sequences of finite 
$R$-modules, then
$R$ \textit{satisfies (ES) for modules}.  
Lemma~\ref{lem:cidim-1} 
implies that one need consider the question in the
case where all three complexes 
have infinite projective dimension.

If $R$ satisfies (ES), then it satisfies (ES) for 
modules;  the converse 
also holds.  In addition, the
rings which satisfy (ES) are exactly those rings for which the converse of 
Corollary~\ref{cor:es} holds.

\begin{thm} \label{thm:es}
For a ring $R$, the following conditions are equivalent.
\begin{enumerate}[(a)]
\item \label{item:es-3}
$R$ satisfies (ES).
\item\label{item:es-4}
 $R$ satisfies (ES) for modules.
\item \label{item:es-5}
Every finite complex of $R$-modules $X$
such that $\cidim_R(X_i)<\infty$ for each integer $i$ satisfies
$\cidim_R(X)<\infty$.
\end{enumerate}
\end{thm}

\begin{proof}
``(\ref{item:es-3})$\implies$(\ref{item:es-5})''.  
Fix a finite complex of $R$-modules $X$
with $\cidim_R(X_i)<\infty$ for each integer $i$.
Since $X$ is finite, proceed by induction on the number
$s$ of modules $X_i$ that are nonzero.   If  $s=0$ or $s=1$, then it 
is immediate that $\cidim_R(X)<\infty$.  If $s>1$, let 
$t=\sup\{i\mid X_i\neq 0\}$ and consider the exact sequence
$0\to X_{\leq t-1}\to X\to \Sigma^tX_t\to 0$.  By induction
$\cidim_R(X_{\leq t-1})<\infty$, 
and since (ES) holds, one has $\cidim_R(X)<\infty$.

``(\ref{item:es-5})$\implies$(\ref{item:es-4})''.  
Let 
$0\to L\stackrel{\nu}{\to} M\stackrel{\phi}{\to} N\to 0$ 
be an exact sequence of nonzero
finite $R$-modules and suppose that two of the modules have 
finite CI-dimension.

Case 1:  $\cidim_R(L),\cidim_R(M)<\infty$.  The complex
$X=0\to L\to M\to 0$
is quasiisomorphic to $N$, and thus, $\cidim_R(N)=\cidim_R(X)<\infty$
by assumption.

Case 2:  $\cidim_R(M),\cidim_R(N)<\infty$.  This is similar to Case 1.

Case 3:  $\cidim_R(L),\cidim_R(N)<\infty$.  Fix a finitely generated 
projective $R$-module $P$ with a surjection 
$\alpha\colon P\twoheadrightarrow N$.
Lemma~\ref{lem:cidim-1} implies 
that $K=\ker(\alpha)$ has $\cidim_R(K)<\infty$.  Let
$\gamma\colon P\to M$ be a map such that $\alpha=\phi\gamma$;
it is straightforward to check that 
there is an exact sequence
\[ 0\to K\to P\oplus L\xrightarrow{(\gamma \,\, \nu)} M\to 0.\]
Lemma~\ref{lem:cidim-1} implies 
that $\cidim_R(P\oplus L)<\infty$.
Since $\cidim_R(K)<\infty$, this implies that 
$\cidim_R(M)<\infty$ by Case 1.

``(\ref{item:es-4})$\implies$(\ref{item:es-3})''.  
Fix an exact sequence 
$0\to X\to Y\to Z\to 0$
of homologically finite complexes of $R$-modules 
such that two of the complexes have finite CI-dimension and all three 
complexes have infinite projective dimension.
By Proposition~\ref{prop:ex-seq-resol}, there exists a commuting diagram
\[ \xymatrix{0\ar[r] & T\ar[r]\ar[d]_{\simeq} & 
U\ar[r]\ar[d]_{\simeq} & V\ar[r]\ar[d]_{\simeq} & 0 \\
0\ar[r] & X\ar[r] & Y\ar[r] & X\ar[r] & 0
} \]
where each row is exact and each vertical map is a 
degreewise finite projective resolution.  
Replace the original sequence with the top row of this diagram 
to assume that each complex is bounded below and 
consists of finitely generated projectives.

For $s=\max\{\sup X, \sup Y, \sup Z\}$ one has an exact sequence
\[ 0\to C_s^X\to C_s^Y\to C_s^Z \to 0.\]
Using Proposition~\ref{prop:cidim-4}, our assumptions 
imply that two of the modules in this sequence have finite CI-dimension.
Since $R$ satisfies property (ES) for modules, the third module also has 
finite CI-dimension.  Using~\ref{prop:cidim-4} again, 
it follows that the third complex in the original sequence has 
finite CI-dimension, as desired.  
\end{proof}

If $R$ satisfies 
(ES) and $X$ is a homologically finite complex whose nonzero
homology modules have finite CI-dimension, then $X$ must also have finite 
CI-dimension.  
Example~\ref{ex:1} below shows that the converse fails.

\begin{prop}
Let $R$ be a ring satisfying (ES) and $X$ a homologically finite complex of 
$R$-modules.  
If $\cidim_R(H_i(X))<\infty$ for all $i$, then $\cidim_R(X)<\infty$.
\end{prop}

\begin{proof}
Since $X$ is homologically finite, argue by induction on 
$s=\sup(X)-\inf(X)$.  
If $s\leq 1$, then $X\simeq\Sigma^jH^j(X)$ for some $j$ and so
$\cidim_R(X)=\cidim_R(H^j(X))+j<\infty$ by~\ref{subprop:cidim-5} 
and~\ref{subprop:cidim-4}.  
When $s>1$, let $t=\sup(X)$ and consider the 
exact sequence 
$0\to Y\to X\to\tau_{\leq t-1}(X)\to 0$.  
By construction, $Y$ and $\tau_{\leq t-1}(X)$ satisfy the induction 
hypothesis and therefore have finite CI-dimension.
As $R$ satisfies (ES), it follows that $\cidim_R(X)<\infty$.
\end{proof}

The following is an example of a ring $R$ and a homologically finite complex 
of $R$-modules $X$ such that $X$ has finite CI-dimension and each 
nonzero homology module $H_i(X)$ has infinite CI-dimension.  
Such a complex must have at least two nonvanishing homology 
modules, and 
this example has exactly two of them.  

\begin{ex} \label{ex:1}
Let $k$ be a field and $R=k[\![S,T]\!]/(S^2,ST,T^2)=k[\![s,t]\!]$ 
with maximal ideal 
$\m=(s,t)R$.  Let $X=(0\to R\stackrel{s}{\to} R\to 0)$.  Then $X$ has 
projective 
dimension 1 and therefore finite CI-dimension.  The homology modules 
are $H_0(X)=R/sR$ and $H_1(X)=\m$.  It is straightforward to verify that 
each of these modules has infinite complexity and therefore cannot have 
finite CI-dimension.
\end{ex}

\section{Lower Complete Intersection Dimension for Complexes} \label{sec:pcidim}

In this section, we consider the lower complete intersection 
dimension, which was introduced for 
modules in~\cite{gerko:ohd} under the name
``polynomial complete intersection dimension''.   We extend this 
dimension to the category of homologically finite complexes and 
present its basic properties.  Most of the results in this 
section have analogues for CI-dimension, and it might seem natural to 
present the two dimensions in the same section.  However, the 
underlying ideas are rather different, so we consider them separately.

We begin with a more general situation coming 
from~\cite[p.~99]{auslander:smt}.  

\begin{defn} \label{def:res-1}
For a ring $R$, a full 
subcategory $\B$ of the category of finite $R$-modules
is a \emph{resolving subclass} if it 
satisfies the following.
\begin{enumerate}[(1)]
\item Every finitely generated projective R-module is in $\B$.
\item \label{item:res-0a}
If $0\to A\to B\to C\to 0$ is an exact sequence of finite $R$-modules
with $C\in\B$, then $A\in\B$ if and only if $B\in\B$.
\item \label{item:res-0b}
If $A,C$ are finite $R$-modules and 
$B=A\oplus C$ is in $\B$, then $A,C\in\B$.
\end{enumerate}
A \emph{$\B$-resolution} of a homologically finite 
complex of $R$-modules
$X$ is a bounded below complex $B\simeq X$ 
with each $B_i$ in $\B$.  
The \emph{$\B$-dimension} of $X$ is
\[ \bdim_R(X)=\inf\{\sup \{ i\mid B_i\neq 0\}\mid \text{$B$
is a $\B$-resolution of $X$}\}.\]
\end{defn}

Certain fact follow from the definition.

\begin{properties} \label{props:res-cat}
Fix a ring $R$ and a homologically 
finite complex of $R$-modules $X$.
\begin{subprops} \label{subprop:res-cat-1}
Each degreewise finite projective resolution of $X$ is a
$\B$-resolution.  
\end{subprops}
\begin{subprops} \label{subprop:res-cat-6}
$\bdim_R(X)\in\{-\infty\}\cup\zz\cup\{\infty\}$.
\end{subprops}
\begin{subprops} \label{subprop:res-cat-2}
$\bdim_R(X)=-\infty$ if and only if $X\simeq 0$.
\end{subprops}
\begin{subprops} \label{subprop:res-cat-4}
Each integer $n$ yields $\bdim_R(\Sigma^nX)=\pcidim_R(X)+n$.
\end{subprops}
\begin{subprops} \label{subprop:res-cat-5}
$\sup(X)\leq \bdim_R(X)$.
\end{subprops}
\end{properties}

\begin{para}
With the previous sections in mind, let $R$ be a ring and set
\[ \mathcal{C}=\{M\mid\cidim_R(M)=0\}\cup\{0\}.\]
One might be tempted to consider the $\mathcal{C}$-dimension arising 
from this choice.  However, in the absence of the property (ES), the 
class $\mathcal{C}$ is not known to be a resolving subclass.  When 
(ES) is satisfied, though, it is straightforward to verify that 
$\mathcal{C}\text{-}\dim_R(X)=\cidim_R(X)$ using 
Propositions~\ref{prop:cidim-4} and~\ref{prop:res-3}.
\end{para}

The following proposition is a version of~\cite[(3.12)]{auslander:smt} 
for complexes.   In the way that Schanuel's lemma allows for the 
computation of $\pd_R(M)$ from an arbitrary projective resolution of a 
module $M$, this result shows that $\bdim_R(X)$ can be computed from 
any $\B$-resolution of $X$.

\begin{prop} \label{prop:res-cat-1}
Let $R$ be a ring and $\B$ a resolving subclass of the category of 
finite $R$-modules.  Consider two complexes of finite
$R$-modules  
\begin{align*}
A&=0\to A_m\to A_{m-1}\to\cdots\to A_n\to 0 \\
B&=0\to B_{m+j}\to B_{m+j-1}\to\cdots\to B_p\to 0 
\end{align*}
with $j\geq 0$ and such that $T\simeq U$ and
$A_{m-1},\ldots\, A_n, B_{m+j-1},\ldots, B_p\in \B$.  
If $A_m$ is in $\B$ then $B_{m+j}$ is in $\B$.
\end{prop}

\begin{proof}
When $\HH(A)=0=\HH(B)$, one uses~\ref{def:res-1}(\ref{item:res-0a})
inductively to show that $\im(\partial_i^A)$ and $\im(\partial_i^B)$
are in $\B$ for $i\leq m$; in particular, both $A_m$ and $B_{m+j}$ are in 
$\B$.

In general, it suffices to consider the case $j=0$.  Indeed, 
since $\sup(B)=\sup(A)\leq m$, one has $\tau_{\geq m}(B)\simeq B$.  By 
the case $j=0$, the module $\tau_{\geq m}(B)_m=C_m^B$ is in $\B$.  
Applying the previous paragraph to the 
exact complex 
\[ 0\to B_{m+j}\to\cdots\to B_{m+1}\to C_m^B\to 0 \]  
one concludes that $B_{m+j}$ is in $\B$.

By Lemma~\ref{lem:back-1}, there exists a degreewise finite free 
resolution
$\sigma\colon P\xrightarrow{\simeq} A$ such that each $\sigma_i$ surjective.  
Since $A\simeq B$, there 
exists a quasiisomorphism $\rho\colon P\xrightarrow{\simeq} B$ 
by~\ref{para:back-3}.
Let $P'=\tau_{\leq m}(P)$ and consider the canonical quasiisomorphism 
$\epsilon\colon P\to P'$.  Because $A_{m+1}=0=B_{m+1}$, it follows that
$\sigma$ and $\rho$ factor through $\epsilon$.  This gives
quasiisomorphisms $\sigma'\colon P'\xrightarrow{\simeq} A$ and
$\rho'\colon P'\xrightarrow{\simeq} B$ such that each
$\sigma_i'$ is surjective.  By construction, the complex 
$P'=0\to P'_m\to P_{m-1}\to\cdots\to P_q\to 0$ has 
$P_{m-1},\ldots,P_q\in\B$.

In order to first see
that $P'_m\in\B$, set $U=\ker(\sigma')$,
which is homologically zero since $\sigma'$ is a quasiisomorphism.  
For $i<m$, applying~\ref{def:res-1}(\ref{item:res-0a}) to
the exact sequence $0\to U_i\to P_i\to A_i\to 0$ implies
that $U_i$ is in $\B$.  Since $\HH(U)=0$ and $U_i=0$ for $i>m$ 
and $i<n$, one has $U_m\in\B$.  The exact sequence 
$0\to U_m\to C_m^P\to A_m\to 0$ implies that $P'_m$ is in $\B$.  

To show that $B_m\in\B$, let $V=\mathrm{Cone}(\rho')$ denote the 
mapping cone of $\rho'$, which is bounded below.  
Since $\rho'$ is a quasiisomorphism, $\HH(V)=0$.  
Since $V_i=B_i\oplus P'_{i-1}$ for each $i$, it follows 
from~\ref{def:res-1}(\ref{item:res-0a}) that
$V_i\in\B$ for $i\leq m-1$.  As above,
one deduces that $\im(\partial_i^V)\in\B$ 
for $i\leq m$.  Furthermore, $V_{m+1}=P'_m$ is in 
$\B$, so the exact sequence
$0\to V_{m+1}\to V_{m}\to\im(\partial_m^V)\to 0$ implies that 
$V_{m}$ is in $\B$.  As $V_m=B_m\oplus P_{m-1}$,
it follows that $B_m$ is in $\B$.  
\end{proof}

One can describe $\bdim(X)$ 
in terms of the inclusion of $C^B_n$ 
in $\B$ for an arbitrary $\B$-resolution $B\simeq X$.

\begin{cor} \label{cor:res-cat-5}
Every $\B$-resolution $B$ of a homologically 
finite complex of $R$-modules $X$ satisfies
\[ \bdim_R(X)=\inf\{n\geq\sup(X)\mid C_n^B\in\B\}.\]
\end{cor}

\begin{proof}
Let $t=\bdim_R(X)$ and $u=\inf\{n\geq\sup(X)\mid C_n^B\in\B\}$.
If $t<\infty$, then $C^B_t$ is in $\B$.  Indeed,
fix a $\B$-resolution $A\simeq X$ with 
$A_i=0$ for all $i>t$.  Then
$\tau_{\leq t}(B)\simeq B\simeq X\simeq A$ since $t\geq\sup(X)$, so
$C^B_t=\tau_{\leq t}(B)_t\in \B$ by Proposition~\ref{prop:res-cat-1}.  

Whether or not $t$ is finite, this shows that 
$t\geq u$.  If $u=\infty$, then $t=u$.
If $u<\infty$,  
then $\tau_{\leq u}(B)$ is a bounded $\B$-resolution of $X$ and 
so $t\leq u$.
\end{proof}

The $\B$-dimension of a complex can be 
computed from that of the syzygies arising from any $\B$-resolution.
Compare this to Proposition~\ref{prop:cidim-4}.

\begin{prop} \label{prop:res-3}
Let $R$ be a ring and $\B$ a resolving subclass of the category of 
finite $R$-modules.  Fix a $\B$-resolution $B$ of a homologically 
finite complex of $R$-modules $X$ and an integer
$n\geq\sup(X)$.
\begin{enumerate}[(i)]
\item
If $C^B_n=0$, then $\bdim_R(X)<n$.
\item
If $C^B_n\neq 0$, then $\bsdim_R(C^B_n)=\max\{0,\bdim_R(X)-n\}$.
\end{enumerate}
\end{prop}

\begin{proof}
Since $n\geq \sup(X)=\sup(B)$, one has 
$\tau_{\leq n}(B)\simeq B\simeq X$.
If $C^B_n=0$, then $\tau_{\leq n}(B)$ is 
a $\B$-resolution of $X$ with $\tau_{\leq n}(B)_i=0$ 
for all $i\geq n$, and it 
follows that $\bdim_R(X)<n$.  Therefore, assume that $C^B_n\neq 0$ and  
let $t=\bdim_R(X)$.  

Case 1:
$t\leq n$.  Corollary~\ref{cor:res-cat-5} implies that
$C^B_t$ is in $\B$.  From~\ref{def:res-1}(\ref{item:res-0a}), it follows that 
$C^B_n$ 
is in $\B$, as well.  Thus, $\bdim_R(C^B_n)=0$ and the formula holds.

Case 2:  $t=\infty$.  From Corollary~\ref{cor:res-cat-5}, it follows that,
for all $m\geq \sup(X)$, the module $C_m^B$ is not in $\B$.  
Since the complex $\Sigma^{-n}(B_{\geq n})$ is a 
$\B$-resolution of $C^B_n$, another application of~\ref{cor:res-cat-5}
yields $\bdim_R(C^B_n)=\infty$, verifying the formula.

Case 3:  $\infty>t>n$.  Again by Corollary~\ref{cor:res-cat-5}, the module 
$C^B_t$ is in $\B$ and $C^B_i$ is not in $\B$ for $i=n,\ldots,t-1$.  
Therefore, the complex
\[ \Sigma^{-n}(\tau_{\leq i}(B))=
0\to C^B_i \to B_{i-1} \to\cdots\to B_n \to 0 \]
is a $\B$-resolution of $C^B_n$ when $i=t$, and is not a 
$\B$-resolution when $i<t$.
By~\ref{cor:res-cat-5}, $\bdim_R(C^B_n)=t-n$ and the formula 
holds.
\end{proof}

\begin{cor} \label{cor:res-3}
For a homologically 
finite complex of $R$-modules $X$, the following conditions are 
equivalent.
\begin{enumerate}[(a)]
\item \label{item:res-cat-4}
$\bdim_R(X)<\infty$.
\item \label{item:res-cat-5}
Each $\B$-resolution $B\simeq X$ and each
$n\geq\sup(X)$ yield $\bdim_R(C_n^B)<\infty$.
\item \label{item:res-cat-6}
Some $\B$-resolution $B\simeq X$ and some 
$n\geq\sup(X)$ yield $\bdim_R(C^B_n)<\infty$. \qed
\end{enumerate}
\end{cor}

The $\B$-dimension behaves well with respect 
to exact sequences of complexes.  As discussed in 
Section~\ref{sec:es}, this is stronger than what we currently know
for CI-dimension.  The corresponding statement for 
\pci-dimension of modules is~\cite[(2.8)]{gerko:ohd}.

\begin{cor} \label{cor:res-cat-4}
An exact sequence of 
homologically finite complexes of $R$-modules
$0\to X^1\to X^2\to X^3\to 0$
and integers $i,j,k$ such that $\{i,j,k\}=\{1,2,3\}$ yield
\[ \bdim_R(X^k)\leq\max\{\bdim_R(X^i),\bdim_R(X^j)\}+1. \]
In particular, if $\bdim_R(X^i)$ and $\bdim_R(X^j)$ are finite, 
then $\bdim_R(X^k)<\infty$.
\end{cor}

\begin{proof}
Almost identical to that of the implication 
``(\ref{item:es-4})$\implies$(\ref{item:es-3})'' in 
Theorem~\ref{thm:es};  use Corollary~\ref{cor:res-3} in place 
of Proposition~\ref{prop:cidim-4}.
\end{proof}

We now specialize the $\B$-dimension to the lower complete 
intersection dimension.
For a nonzero finite module, considered as a complex concentrated in 
degree 0, the definition is the same as
that given in~\cite[(2.3)]{gerko:ohd}.

\begin{defn} \label{def:pcidim-1}
Let $R$ be a ring.  The \emph{\pci-class of $R$}, denoted 
$\pciclass(R)$, is the collection of 
totally reflexive $R$-modules $T$ such that,
for every maximal ideal $\m$ of $R$, the localized
module $T_{\m}$ has finite complexity over $R_{\m}$.  
Thus, a finite module $T$ is in $\pciclass(R)$ if and only if, 
for every maximal ideal $\m$ of $R$, the $R_{\m}$-module $T_{\m}$ is 
totally reflexive and has finite complexity.

From~\cite[(1.1.10,11)]{christensen:gd} and~\cite[(4.2.4)]{avramov:ifr} 
it follows that $\pciclass(R)$ is a resolving subclass of the category 
of finite $R$-modules.  The resulting homological dimension
$\pcidim_R$ is the \emph{lower complete intersection dimension}.
\end{defn}

Of course, the results stated for $\B$-dimension hold for 
\pci-dimension.  We continue with properties specifically
for the \pci-dimension.  The first of these states that like 
CI-dimension (\ref{prop:cidim-2}) 
the  \pci-dimension of a complex does not increase after
localizing and is determined locally.  The result for 
finite modules is~\cite[(2.11)]{gerko:ohd}.

\begin{prop} \label{prop:pcidim-1}
Let $R$ be a ring and $X$ a homologically finite complex 
of $R$-modules.  For every multiplicative subset
$S\subset R$ there is an inequality
\[ \pcidim_{S^{-1}R}(S^{-1}X)\leq\pcidim_R(X).\]
Furthermore, there are equalities
\begin{align*}
\pcidim_R(X)&=\sup\{\pcidim_{R_{\m}}(X_{\m})\mid \m\in\mspec(R)\} \\
&=\sup\{\pcidim_{R_{\p}}(X_{\p})\mid \p\in\spec(R)\}.
\end{align*}
\end{prop}

\begin{proof}
The inequality follows readily;  use~\cite[(2.11)]{gerko:ohd} to 
show that a \pci-resolution of $X$ over $R$ localizes to a 
\pci-resolution of $S^{-1}X$ over $S^{-1}R$.

For the other formulas, set
$v=\sup\{\pcidim_{R_{\p}}(X_{\p})\mid \p\in\spec(R)\}$.
It follows from the inequality above that we need only verify that
$\pcidim_R(X)\leq v$.  To this end,
assume that $v<\infty$.  
Fix a \pci-resolution $U\simeq X$ over $R$ and note that
\[ v\geq\sup\{\sup(X_{\p})\mid\p\in\spec(R)\}=\sup(X).\]
For every $\p$, 
the complex $U_{\p}$ is a \pci-resolution of $X_{\p}$ over 
$R_{\p}$ and $C_u^{U_{\p}}\cong (C_u^U)_{\p}$.  By 
Corollary~\ref{cor:res-cat-5}, the module $(C_u^U)_{\p}$ is in 
$\bs(R_{\p})$ 
for all $\p$.  By definition, $C_u^U\in\bs(R)$, so that 
$\pcidim_R(X)\leq u$.
\end{proof}

The following result explains the position of \pci-dimension in 
the hierarchy of homological dimensions and shows that complexes of 
finite \pci-dimension over a local ring satisfy an AB-formula.  
That this holds for finite modules is in~\cite[(2.6,7)]{gerko:ohd}.  
It is important to note that each of the given inequalities can be 
strict.  For the first and third inequalities, this is 
straightforward.  For the second inequality, 
this is due to Veliche~\cite[Main Theorem (4)]{veliche:cmfhd}.

\begin{prop} \label{prop:pcidim-3}
Let $R$ be a ring and $X$ a homologically finite complex of 
$R$-modules.  There are inequalities
\[ \gdim_R(X)\leq\pcidim_R(X)\leq\cidim_R(X)\leq\pdim_R(X); \]
when one of these dimensions is finite it is equal to those 
on its left.    If $R$ is 
local and $\pcidim_R(X)<\infty$, then 
$\pcidim_R(X)=\depth(R)-\depth_R(X)$.
\end{prop}

\begin{proof}
By Proposition~\ref{prop:pcidim-1},
it suffices to consider the case when $R$ is local.
The third inequality is in Proposition~\ref{prop:cidim-1}.

The first inequality holds because every \pci-resolution of $X$ is a 
G-resolution.  
When $\pcidim_R(X)<\infty$, let $T\simeq X$ be a \pci-resolution.
For every $n\geq\sup(X)$, one has $\pcidim_R(C^T_n)<\infty$ 
by Proposition~\ref{prop:res-3}.  The AB-formulas~\ref{para:AB}
and~\cite[(2.7)]{gerko:ohd} imply the equality
$\gdim_R(C^T_n)=\pcidim_R(C^T_n)$
and it follows that $C^T_n$ is in $\pciclass(R)$ if and only if it is 
totally reflexive.  Corollary~\ref{cor:res-cat-5} and the corresponding 
equality for G-dimension~\cite[(2.3.7)]{christensen:gd}, 
\[ \gdim_R(X)=\inf\{n\geq\sup(X)\mid\text{$C^T_n$ is totally 
reflexive}\} \]
imply that 
$\pcidim_R(X)=\gdim_R(X)$.  From the AB-formula~\ref{para:AB}
it follows that this equals $\depth(R)-\depth_R(X)$.

For the second inequality, assume that
$\cidim_R(X)<\infty$.  Using the AB-formula, it suffices to show 
that $\pcidim_R(X)<\infty$.
Let $F\simeq X$ be a degreewise finite free resolution.  
By~\ref{prop:cidim-4}, 
one has $\cidim_R(C^F_q)\leq 0$ for $q\gg 0$.  
Thus, $\pcidim_R(C^F_q)\leq 0$ by~\cite[(2.6)]{gerko:ohd}, i.e., 
$C^F_q\in\bs(R)$, 
and~\ref{cor:res-cat-5} implies that
$\pcidim_R(X)<\infty$.
\end{proof}

The next result is the analogue of 
Proposition~\ref{prop:cidim-3} for \pci-dimension.  The local case 
for modules is given in~\cite[(2.5)]{gerko:ohd}.

\begin{prop} \label{prop:pcidim-4}
For a ring $R$ with $\dim(R)<\infty$
the following are equivalent.
\begin{enumerate}[(a)]
\item \label{item:pcidim-1}
$R$ is locally a complete intersection.
\item \label{item:pcidim-2}
Each homologically finite complex of $R$-modules $X$ satisfies
\[ \pcidim_R(X)\leq\dim(R)+\sup(X).\]
\item  \label{item:pcidim-3}
Each maximal ideal $\m\subset R$ satisfies
$\pcidim_R(R/\m)<\infty$.
\end{enumerate}
\end{prop}

\begin{proof}
``(\ref{item:pcidim-1})$\implies$(\ref{item:pcidim-2})''.  
For a homologically finite complex of $R$-modules $X$, one has
\[ \pcidim_R(X)\leq\cidim_R(X)\leq\dim(R)+\sup(X)<\infty \]
where the first inequality is by Proposition~\ref{prop:pcidim-3} and 
the second is by Proposition~\ref{prop:cidim-3}.

``(\ref{item:pcidim-2})$\implies$(\ref{item:pcidim-3})'' is trivial.

``(\ref{item:pcidim-3})$\implies$(\ref{item:pcidim-1})''.  
One has
$\pcidim_{R_{\m}}(R_{\m}/\m R_{\m})=\pcidim_R(R/\m)$ for each $\m$ by 
Proposition~\ref{prop:pcidim-1}, and so  
$R_{\m}$ is a complete intersection by~\cite[(2.5)]{gerko:ohd}.
\end{proof}

The complexes of finite \pci-dimension 
are exactly those that behave as a whole like the modules in the 
\pci-class.

\begin{thm} \label{thm:pcidim-1}
A homologically finite complex $X$ over a ring $R$ has finite
\pci-dimension if and only if 
$\gdim_R(X)$ is finite and $\cx_{R_{\m}}(X_{\m})$ is finite for all maximal 
ideals $\m$ of $R$.
\end{thm}

\begin{proof}
Let $P\simeq X$ be a degreewise finite $R$-projective resolution.   

Assume first that $p=\pcidim_R(X)<\infty$.  Then $\gdim_R(X)<\infty$ by 
Proposition~\ref{prop:pcidim-3}.  
The module $C_p^{P_{\m}}\cong (C_p^P)_{\m}$ is in $\pciclass(R_{\m})$ by 
Proposition~\ref{prop:pcidim-1} and Corollary~\ref{cor:res-cat-5}.  
The result now follows because
$\cx_{R_{\m}}(X_{\m})=\cx_{R_{\m}}((C_p^P)_{\m})<\infty$.

Assume now that $g=\gdim_R(X)<\infty$ 
and $\cx_{R_{\m}}(X_{\m})<\infty$ for all maximal 
ideals $\m$ of $R$.  The module $C_g^P$ is totally reflexive over $R$ 
by~\cite[(2.3.7)]{christensen:gd}.  For all $\m$, one has
$\cx_{R_{\m}}((C_g^P)_{\m})=\cx_{R_{\m}}(X_{\m})<\infty$.  Hence,
$C_g^P$ is in $\pciclass(R)$ and it follows that 
$\pcidim_R(X)<\infty$.
\end{proof}

A souped-up version of
Corollary~\ref{cor:cidim-2}(\ref{item:cidim-7}) 
is satisfied by \pci-dimension.

\begin{prop} \label{prop:pcidim-5}
Let $\vf\colon R\to S$ be a flat ring homomorphism
and $X$ a homologically finite complex of $R$-modules.  There is an 
inequality
\[ \pcidim_S(X\otimes_R S)\leq\pcidim_R(X) \]
with equality when $\vf$ is faithfully flat.
\end{prop}

\begin{proof}
For any $M\in\bs(R)$, it follows from flatness that
$M\otimes_R S$ is in $\bs(S)$.  Thus, a 
\pci-resolution of $X$ over $R$ base-changes to a \pci-resolution of 
$X\otimes_R S$ over $S$, and hence the inequality holds.

When $\vf$ is faithfully flat and $M$ is a finite $R$-module, it 
follows readily that $M$ is in $\pciclass(R)$ if and only if  
$M\otimes_R S$ is in $\pciclass(S)$.  
To show that $\pcidim_S(X\otimes_R S)=\pcidim_R(X)$, fix
a \pci-resolution $U\simeq X$ over $R$.  
Then $U\otimes_R S$ is a \pci-resolution of $X\otimes_R S$ over $S$, and 
$C_n^{U\otimes_R S}\cong C_n^U\otimes_R S$ for each integer $n$.
Furthermore, $\sup(X\otimes_R S)=\sup(X)$, so one has
\begin{align*}
\pcidim_S(X\otimes_R S)
 &=\inf\{n\geq\sup(X\otimes_R S)\mid C_n^U\otimes_R S\in\bs(S)\}\\
 &=\inf\{n\geq\sup(X)\mid C_n^U\in\bs(R)\}\\
 &=\pcidim_R(X)
\end{align*}
where the first and third equalities are by
Corollary~\ref{cor:res-cat-5}.
\end{proof}

The following is a version of 
Corollary~\ref{cor:cidim-2}(\ref{item:cidim-8}) for \pci-dimension.

\begin{prop} \label{prop:pcidim-6}
Let $Q\to R$ be a surjective ring homomorphism with kernel generated 
by a $Q$-regular sequence of length $c$.  
Every homologically finite complex of 
$R$-modules $X$ satisfies
\[ \pcidim_Q(X)=\pcidim_R(X)+c.\]
In particular, $\pcidim_Q(X)$ is finite if and only if 
$\pcidim_R(X)$ is finite.
\end{prop}

\begin{proof}
By Proposition~\ref{prop:pcidim-1}, it suffices to consider the case 
where $Q$ and $R$ are local.  By~\cite[(2.3.12)]{christensen:gd},
\ref{prop:pcidim-3}, and~\ref{thm:pcidim-1}, 
one needs only show that 
$\cx_R(X)$ and $\cx_Q(X)$ are simultaneously finite.  Assume that 
$\HH(X)$ is nonzero and fix 
a degreewise finite free resolution $P\simeq X$ and an integer 
$n\geq\sup(X)$.
The complex $P_{\leq n-1}$ has finite projective dimension 
over $R$, and thus also over $Q$.  The exact sequence
$0\to P_{\geq n}\to P\to P_{\leq n-1}\to 0$ implies that
\[ \cx_Q(X)=\cx_Q(P)=\cx_Q(P_{\geq n})=\cx_Q(C^P_n) \]
and 
similarly, $\cx_R(X)=\cx_R(C^P_n)$.  Thus, it suffices to consider the 
case where $X$ is a module.  This case is in~\cite[(5.2.4)]{avramov:cid}.
\end{proof}

The final result of this section is the analogue of 
Theorem~\ref{thm:cidim-stability} for \pci-dimension.

\begin{thm} \label{thm:pcidim-stability}
Let $R$ be a ring and $X,P$ homologically finite complexes of 
$R$-modules.  
If $\pd_R(P)$ is finite, 
then
\begin{align*}
\pcidim_R(X\lotimes_R P)&=\pcidim_R(X)+\pcidim_R(P) & \text{and} \\
\pcidim_R(\rhom_R(P,X))&=\pcidim_R(X)-\inf(P).
\end{align*}
In particular, the \pci-dimensions of the complexes
$X$, $X\lotimes_R P$, and $\rhom_R(P,X)$
are simultaneously finite.
\end{thm}

\begin{proof}
As in the proof of 
Theorem~\ref{thm:cidim-stability}, it suffices to show that 
the complexes $X$, $X\lotimes_R P$, and $\rhom_R(P,X)$ have finite
\pci-dimensions simultaneously
when $R$ 
is local. 
By~\cite[(5.1,7)]{iyengar:gdolh}, 
the G-dimensions of the complexes
$X$, $X\lotimes_R P$, and $\rhom_R(P,X)$
are simultaneously finite,
so it suffices to show that 
\[ \cx_R(X\lotimes_R P)=\cx_R(X)=\cx_R(\rhom_R(P,X)).\]
The first equality follows from Lemma~\ref{para:back-cx-add}.  
This lemma also implies the second equality because of the
isomorphism $\rhom_R(P,X)\simeq X\lotimes_R\rhom_R(P,R)$ and
since $\rhom_R(P,R)$ is homologically finite and
$\pd_R(\rhom_R(P,R))$ is finite.
\end{proof}

\section{Global Homological Dimensions} \label{sec:gldim}

We use the homological 
dimensions discussed in the previous sections to define global 
homological dimensions of rings similar to the 
global dimension of~\cite{cartan:ha}.  
The primary focus is the CI-dimension.
The first proposition of this section motivates our definition of the 
global CI-dimension of a ring $R$.  Similar results hold for 
\pci-dimension and G-dimension.

\begin{prop} \label{prop:ghd-1}
For a ring $R$ and an integer $n$, the following are 
equivalent.
\begin{enumerate}[(a)]
\item \label{item:ghd-1}
Each homologically finite complex of $R$-modules $X$ satisfies
\[ \cidim_R(X)\leq n+\sup(X). \]
\item \label{item:ghd-2}
Each finite $R$-module $M$ satisfies $\cidim_R(M)\leq n$.
\end{enumerate}
\end{prop}

\begin{proof}
The implication ``(\ref{item:ghd-1})$\implies$(\ref{item:ghd-2})'' is 
clear.  
For the other implication, assume (\ref{item:ghd-2}) holds and 
fix a homologically finite complex of $R$-modules $X$.  Set
$s=\sup(X)$, and let $P\simeq X$ be a degreewise finite projective 
resolution.  Then $\cidim_R(C^P_s)\leq n$, by assumption, and 
Proposition~\ref{prop:cidim-4} 
implies that $\cidim_R(X)-s\leq n$.  
\end{proof}

\begin{defn} \label{def:ghd-1}
For a ring $R$, the \emph{global CI-dimension of $R$} is 
\[ \glcidim(R):=\inf\{n\in\zz\mid\text{$\cidim_R(M)\leq n$, $\forall$ 
finite $R$-modules $M$}\}. \]
The above proposition implies that this is equal to
\[ \inf\{ n\in\zz\mid\text{$\cidim_R(X)\leq n+\sup(X)$, $\forall$ 
homologically finite $R$-complexes $X$}\}.\]
\end{defn}

In a similar manner, 
one can define the global \pci-dimension and global 
G-dimension.  Each of these quantities is 
in $\nn\cup\{\infty\}$.

The hierarchy of global homological dimensions follows from 
Proposition~\ref{prop:pcidim-3}.

\begin{prop} \label{prop:ghd-2}
For a ring $R$, there are inequalities
\[ \glgdim(R)\leq\glpcidim(R)\leq\glcidim(R)\leq\gldim(R); \]
when one of these dimensions is finite it is equal to those 
on its left.  \qed
\end{prop}

Like the CI-dimension, the global CI-dimension is 
determined locally.

\begin{prop} \label{prop:ghd-3}
For a ring $R$, there are (in)equalities
\begin{align*}
\dim(R)
  &\leq\sup\{\cidim_R(R/\m)\mid\m\in\mspec(R)\}\\
  &=\sup\{\glcidim(R_{\m})\mid\m\in\mspec(R)\} \\
  &= \glcidim(R)
\end{align*}
with equality in the first spot when $\glcidim(R)<\infty$.
\end{prop}

\begin{proof}
Set
\begin{align*}
u&=\sup\{\cidim_R(R/\m)\mid\m\in\mspec(R)\}\\
v&=\sup\{\glcidim(R_{\m})\mid\m\in\mspec(R)\} \\
w&=\glcidim(R).
\end{align*}
To verify the inequality $\dim(R)\leq u$, assume that
$u$ is finite.  For each maximal ideal $\m$, one has
$\cidim_{R_{\m}}(R_{\m}/\m R_{\m})=\cidim_R(R/\m)<\infty$.  
By Proposition~\ref{prop:cidim-3}, each $R_{\m}$ is a complete 
intersection, and it follows that
\begin{align*}
\dim(R)
  &=\sup\{\dim(R_{\m})\mid\m\in\mspec(R)\}\\
  &=\sup\{\depth(R_{\m})\mid\m\in\mspec(R)\}\\
  &=\sup\{\cidim_{R_{\m}}(R_{\m}/\m R_{\m})\mid\m\in\mspec(R)\} \\
  &=u.
\end{align*}

Next, we verify the inequalities $u\leq v\leq w\leq u$.  
That $u\leq v$ comes from the inequality
$\cidim_{R_{\m}}(R_{\m}/\m R_{\m})\leq\glcidim(R_{\m})$.
That $v\leq w$
is also straightforward: every finite $R_{\m}$-module is 
of the form $M_{\m}$ for some finite $R$-module $M$ and 
$\cidim_{R_{\m}}(M_{\m})\leq\cidim_R(M)$, so 
$\glcidim(R_{\m})\leq\glcidim(R)$.  

For the final inequality, assume that $u<\infty$.  Then
$R$ is locally a complete intersection, as above.
When $M$ is a 
finite $R$-module, one has
\begin{align*}
\cidim_R(M)
  &=\sup\{\cidim_{R_{\m}}(M_{\m})\mid\m\in\mspec(R)\}\\
  &=\sup\{\depth(R_{\m})-\depth_{R_{\m}}(M_{\m})\mid\m\in\mspec(R)\}\\
  &\leq\sup\{\depth(R_{\m})\mid\m\in\mspec(R)\}\\
  &=\sup\{\cidim_{R_{\m}}(R_{\m}/\m R_{\m})\mid\m\in\mspec(R)\}\\
  &=w.
\end{align*}
By definition, it follows that $\glcidim(R)\leq w$.
\end{proof}

In the same way that the regular rings are characterized
as the rings of finite global dimension, the local complete 
intersection rings of finite Krull dimension are exactly the rings of 
finite global complete intersection dimension.

\begin{thm} \label{thm:ghd-1}
For a ring $R$, the following conditions are equivalent:
\begin{enumerate}[(a)]
\item  \label{item:ghd-3}
$\glcidim(R)=\dim(R)<\infty$;
\item  \label{item:ghd-4}
$\glcidim(R)<\infty$;
\item  \label{item:ghd-5}
$\glpcidim(R)<\infty$;
\item  \label{item:ghd-6}
$R$ is locally a complete intersection and $\dim(R)<\infty$.
\end{enumerate}
\end{thm}

\begin{proof}  
``(\ref{item:ghd-3})$\implies$(\ref{item:ghd-4})'' is trivial.
``(\ref{item:ghd-4})$\implies$(\ref{item:ghd-5})'' follows from 
Proposition~\ref{prop:ghd-2}.

``(\ref{item:ghd-5})$\implies$(\ref{item:ghd-6})''.
Since $\glpcidim(R_{\m})\leq\glpcidim(R)<\infty$, 
Proposition~\ref{prop:pcidim-4} implies that 
$R$ is locally a complete intersection.  Arguing 
as in Proposition~\ref{prop:ghd-3} one sees that
$\dim(R)\leq\glpcidim(R)<\infty$.

``(\ref{item:ghd-6})$\implies$(\ref{item:ghd-3})''  
Proposition~\ref{prop:cidim-3} implies that
$\glcidim(R)\leq \dim(R)<\infty$.  By Proposition~\ref{prop:ghd-3}, 
$\glcidim(R)=\dim(R)$.
\end{proof}

\begin{cor} \label{cor:ghd-1}
Every ring $R$ satisfies $\glpcidim(R)=\glcidim(R)$.
\qed
\end{cor}

Nagata~\cite[A1. Example 1]{nagata:lr} constructed a ring that is 
locally regular with infinite global dimension.  This shows 
that the implication ``locally CI$\implies\glcidim(R)<\infty$'' 
does not hold 
without the additional hypothesis ``$\dim(R)<\infty$''.

The final result of this paper is a version of Theorem~\ref{thm:ghd-1} for 
G-dimension.

\begin{thm} \label{thm:ghd-2}
For a ring $R$, the following conditions are equivalent:
\begin{enumerate}[(a)]
\item  \label{item:ghd-6a}
$\glgdim(R)=\dim(R)<\infty$;
\item  \label{item:ghd-7}
$\glgdim(R)<\infty$;
\item  \label{item:ghd-8}
$R$ is locally Gorenstein and $\dim(R)<\infty$.
\item  \label{item:ghd-9}
$\id_R(R)=\dim(R)<\infty$.
\item  \label{item:ghd-10}
$\id_R(R)<\infty$.
\end{enumerate}
\end{thm}

\begin{proof}
The equivalence of (\ref{item:ghd-6a}), (\ref{item:ghd-7}), and 
(\ref{item:ghd-8}) is verified as in Theorem~\ref{thm:ghd-1}.  The
implication ``(\ref{item:ghd-9})$\implies$(\ref{item:ghd-10})'' is
trivial.  For 
the other equivalences, recall the 
following fact~\cite[(3.5)]{roberts:hiomocr}:  If 
$I$ is a minimal $R$-injective resolution for $R$ and $\m$ is a maximal 
ideal of $R$, the the localized complex $I_{\m}$ is a minimal injective
resolution of $R_{\m}$.  

``(\ref{item:ghd-8})$\implies$(\ref{item:ghd-9})''.  Let $I$ be a 
minimal injective $R$-resolution of $R$;  then 
\begin{align*}
\id_R(R)&=\sup(I) \\
 & = \sup\{\sup(I_{\m})\mid\m\in\mspec(R)\} \\
 & = \sup\{\dim(R_{\m})\mid\m\in\mspec(R)\} \\
 & =\dim(R).
\end{align*}

``(\ref{item:ghd-10})$\implies$(\ref{item:ghd-8})''.  
The chain of inequalities 
\[ \dim(R_{\m})\leq\id_{R_{\m}}(R_{\m})\leq\id_R(R)<\infty \] 
implies that $R$ is locally Gorenstein and 
$\dim(R)\leq\id(R)<\infty$.
\end{proof}

\providecommand{\bysame}{\leavevmode\hbox to3em{\hrulefill}\thinspace}

\end{document}